\documentstyle[12pt]{article}
\author{R.R.Zapatrine}
\date{} 

\title{Les espaces duaux pour les ensembles ordonn\'es arbitraires} 

\def\dff{\sc} 
\def\C{{\bf C}} 
\def\mor{{\tt Mor}} 
\def\b{{\bf B}_2} 
\def\clos{{\tt clos}} 
\def\endproof{\hspace*{\fill}$\Box$}
\def\fix{{\tt Fix}}
\newtheorem{th}{Th\'eor\`eme} 
\newtheorem{lemma}[th]{Lemma} 

\begin{document}
\maketitle 

\begin{abstract} 
On \' etablit une dualit\' e entre les ensembles ordonn\' es et 
certains espaces \`a deux fermetures. 
\end{abstract} 

\section{Introduction} 

En 1982 R. Mayet a \'etabli la dualit\'e entre les 
ensembles ordonn\'es orthocompl\'ement\'es et certains espaces de 
fermeture (dites C-espaces). Suivant \cite{mayet} tout ensemble 
ordonn\'e orthocompl\'ement\'e peut \^etre represent\'e comme 
l'ensembles des parties \`a la fois ferm\'ees et ouvertes. 

Dans cet article cette dualit\'e est \'etendu aux ensembles 
ordonn\'es arbitraires. Le r\'esultat a la forme suivante. Avec 
tout ensemble ordonn\'e $P$ on associe un espace \`a deux 
fermetures $X$ tel que $P$ est isomorphique \`a la collection des 
parties de $X$ ferm\'ees par rapport d'une fermeture et ouvertes 
par rapport d'autre fermeture. 

L'id\'ee principale est d'utiliser la procedure de duplication de 
l'ensemble ordonn\'e $P$. On construit la somme horizontale de $P$ 
et son inverse $P^{op}$ (\ref{eee}), la munit avec la 
compl\'ementation \`a le fa\c con canonique, et alors utilise les 
r\'esultats de \cite{mayet}.

\section{Espaces de fermeture} 

Un {\dff espace de fermeture}  est un couple $(X,\C)$, X \'etant  
un ensemble et $\C:\exp X\to \exp X$ une fermeture sur $X$. Les 
parties de $X$ de la forme $A=\C A$ sont dites 
{\dff $\C$-ferm\'es} et leurs compl\'emants sont dites 
{\dff $\C$-ouverts}. D\'esignons 
\[
C(X)=\{A\subseteq X\mid\quad A\hbox{ est }\C\hbox{-ferm\'e}\} 
\] 
\[ 
O(X)=\{B\subseteq X\mid\quad B\hbox{ est }\C\hbox{-ouvert}\} 
\]
\[ CO(X) = C(X)\cap O(X) 
\]

Une collection $K\subseteq\exp X$ des parties de $X$ est dite la 
{\dff base de la fermeture}\ $\C$, not\'e  
\[ 
\C = \clos(K) 
\] 
si tout $\C$-ferm\'e $A$ est une intersection des \'el\'ements de $K$:  
\[ 
A\in C(X) \Leftrightarrow \exists I \quad A=\cap_{i\in I}K_i 
\hbox{ , }K_i\in K
\] 

Pour deux espaces \`a femeture $(X,\C)$, $(X',\C')$, une 
application $f:X\to X'$ est dite {\dff continue} si pour tout 
$\C'$-ferm\'e $A'$ son preimage $f^{-1}(A')$ est 
$\C$-ouvert-ferm\'e.  Deux espaces \`a fermeture $(X,\C)$ et 
$(X',\C')$ sont dites {\dff hom\'eomorphiqiue} s'il existe un 
couple des bijections continues $X\to X'$ et $X'\to X$. 

Un {\dff es\-p\-a\-ce \`a deux fer\-me\-tu\-res} est un tri\-ple 
$(X,\C_1,\C_2)$ tel que $(X,\C_1)$ aus\-si que $(X,\C_2)$ sont des 
es\-pa\-ces de fer\-me\-tu\-re. D\'e\-s\-i\-gn\-ons pour $i,j=1,2$:  

\[ C_i(X)=\{A\subseteq X\mid\quad A\hbox{ est 
}\C_i\hbox{-ferm\'e}\} \] 
\[ O_j(X)=\{B\subseteq X\mid\quad B\hbox{ 
est }\C_i\hbox{-ouvert}\} \] 
\[C_iO_j(X) = C_i(X)\cap O_j(X) 
\] les 
es\-pa\-ces \`a deux fer\-me\-tu\-res $(X,\C_1,\C_2)$ et 
$(X',\C_1',\C_2')$ sont dites {\dff ho\-m\'e\-o\-mor\-phi\-que} si 
les couples des es\-pa\-ces de fer\-me\-tu\-res $(X,\C_1)$, 
$(X',\C_1')$ aus\-si que $(X,\C_2)$, $(X',\C_2')$ sont 
ho\-m\'e\-o\-mor\-phi\-ques. 

\section{Ensembles ordonn\'es} 

Un ensemble ordonn\'e (EO) est un ensemble non-vide $P$ muni d'une 
relation d'ordre. Un EO $P$ est dit {\dff born\'e} (EOB) s'il 
poss\`ede le plus grand \'el\'ement $1$ est le plus petit 
\'el\'ement $0$. Un EOB $E$ est dit {\dff ensemble ordonn\'e 
compl\'ement\'e} (EOC) s'il est muni d'une loi unaire $p\to 
p^\perp$ telle que, quels que soient $p,q\in E$ 
\begin{equation} 
p^{\perp\perp}=p \label{invol} 
\end{equation} 
\begin{equation} 
p\lor p^\perp=1\qquad p\land 
p^\perp=0 \label{compl} 
\end{equation} 
\begin{equation} 
p\le q\quad\hbox{ entra\^\i ne }\quad 
q^\perp\le p^\perp \label{antiton} 
\end{equation} 

Soient $P,Q$ deux EO. Une application $f:P\to Q$ est dite {\dff 
croissante} si pour tous $p,q\in P$ \ $p\le q$ implique $f(p)\le 
f(q)$, et {\dff anti-croissante} si $p\le q$ implique $f(q)\le 
f(p)$. On designe 
\begin{equation} 
\mor_{EO}(P,Q) = \{f:P\to Q\mid\quad f\hbox{ est croissante}\} 
\label{moreo} 
\end{equation} 

Pour tout EO $P$ on d\'efinit son {\dff inverse} $P^{op}$ \'etant 
une copie de l'ensemble $P$ avec l'ordre inverse. La bijection 
naturelle 
\begin{equation} 
\rho:P\to P^{op} \label{rho} 
\end{equation} 
est anti-croissante. 

Pour deux EOC $E,F$, une application croissante $f:E\to F$ est dite 
{\dff orthocroissante} si $f$ pr\'eserve les compl\'ements: 
$f(p^\perp)=(f(p))^\perp$. D\'esignons 
\begin{equation} 
\mor_{EOC}(E,F)=\{f:E\to F\mid\quad f\hbox{ est orthocroissante}\} 
\label{moreoc}  
\end{equation} 

On re\-ma\-r\-que que pour tout es\-pa\-ce \`a deux 
fer\-me\-tu\-res $(X,\C_1,\C_2)$ l'en\-sem\-b\-le $C_1O_2(X)$ est 
un EO, et pour tout es\-pa\-ce de fer\-me\-tu\-re $(X,\C)$ 
l'en\-sem\-b\-le $CO(X)$ est un EOC par rap\-port de 
com\-p\-l\'e\-men\-ta\-tion des parties de $X$. 

\section{Espaces duaux de Mayet} 

Soit $\b$ un EOC \`a deux \'el\'ements: $\b=\{0,1\}$, $1^\perp=0, 
0^\perp=1$. Soit $E$ un EOC, et introduisons son espace dual 
\begin{equation} 
X=\mor_{EOC}(P,\b) \label{xxx} 
\end{equation} 
Pour tout $p\in E$ posons 
\[ \sigma(p)=\{x\in X\mid\quad x(p)=1\} \] 
et munissons l'ensemble $X$ avec la fermeture $\C$ engendr\'ee par 
$\sigma(E)$ 
\begin{equation} 
\C=\clos(\sigma(E)) \label{ce} 
\end{equation} 

\begin{th} \label{thmayet} Soit $E$ un EOC et $X$ son espace dual 
(\ref{xxx}).  Alors \begin{enumerate} \item $(X,\C)$ et le C-espace 
\item Les EOC $E$ est $CO(X)$ sont isomorphique. Cet isomorphisme 
est r\'ealis\'e par l'application $\sigma: E\to \exp X$ et son 
inverse 
\end{enumerate} 
\end{th} 
\paragraph{Epreuve.} Voir \cite{mayet}. \endproof 
\paragraph{Corollaire.} Les parties $CO(X)$ ne sont que les images 
$\sigma(p)$ des \'el\'ements $p\in E$. 

\section{La dualit\'e pour les ensembles ordonn\'es arbitraires} 

Maintenant soit $P$ un ensemble ordonn\'e arbitraire, et soit 
$Q=P^{op}$ son inverse avec l'anti-iso\-mor\-phisme canonique 
$\rho$ (\ref{rho}).  Formons l'EOC engendr\'e par $P$:  
\begin{equation} 
E=P\cup Q\cup\{0,1\} \label{eee} 
\end{equation} 
avec l'ordre suivant. 1 est le plus grant \'el\'ement de $E$, 
0 est le plus petit, et pour tous $p,p'\in P\quad  p\le p'$ dans 
$E$ si et seulement si $p\le p'$ dans $P$, pour tous $q,q'\in 
Q\quad  q\le q'$ dans $E$ si et seulement si 
$\rho^{-1}(q)\le_P\rho^{-1}(q')$ dans $P$ et tout couple 
$p\in P, q\in Q$ n'est pas comparable. Munissons l'espace 
$X=\mor_{EOC}(E,B_2)$ (\ref{moreoc}) avec la fermeture $\C$ 
(\ref{ce}).  D\'esignons 
\[ Y=\mor_{EO}(P,B_2) \] 
Soient pour tout $p\in P$: 
\begin{equation} 
\begin{array}{rcl}
\sigma_1(p)&=&\{y\in Y\mid\quad y(p)=1\}\cr 
\sigma_2(p)&=&\{y\in Y\mid\quad y(p)=0\} 
\end{array} \label{s12} 
\end{equation} 
et munissons l'espace $Y$ avec un couple 
des fermetures $\C_1,\C_2$ et avec la fermeture $\C$: 
\begin{equation} 
\begin{array}{rcl}
\C_1&=&\clos\{\sigma_1(P)\}\cr 
\C_2&=&\clos\{\sigma_2(P)\}\cr 
\C=&\C_1\lor\C_2=&\clos\{\sigma_2(P)\cup\sigma_2(P)\}
\end{array} \label{c12} 
\end{equation} 

\begin{lemma} Les espaces de fermeture $X,\C$ et $Y,\C$ sont 
hom\-\'e\-o\-mor\-phi\-ques. 
\end{lemma} 
\paragraph{Epreuve.} Premi\`erement \'etablissons la bijection 
entre les ensembles $X$ et $Y$. Pour tout $x\in X$ sa restriction 
$x^P=x\mid_P:P\to B_2$ est un \'element de $Y$. De contraire, 
soit $y\in Y$. Posons $y^E(0)=0$ et $y^E(1)=1$, pour tout $p\in 
P\subseteq E$ soit $y^E(p)=y(p)$ et pour tout $q\in Q$ (\ref{eee}) 
posons $y^E(q)=1-y(\rho^{-1}(q))$, alors cet expansion $y^E$ de 
$y$ devient l'\'el\'ement de $X$ \'etant orthocroissant. 

Pour \'etablir que les bijections $x\mapsto x^P$ et 
$y\mapsto y^E$ sont bicontinues il faut v\'erifier que les 
pr\'eimages des \'el\'ements des bases des fermetures 
appropri\'ees sont ferm\'es. Commen\c cons avec $\phi:x\mapsto 
x^P:X\to Y$.  Les \'el\'ements de la base de $\C(Y)$ sont les 
parties de $Y$ de deux types (\ref{s12}). Pout tout $\sigma_1(p)$ 
nous avons \( \phi^{-1}(\sigma_1(p))=\sigma(p)\in C(X) \), pour 
tout $\sigma_2(p)$ nous avons \( \phi^{-1}(\sigma_2(p))=\{x\in 
X\mid\quad x(p)=0\}=\sigma(p^\perp)\in C(X) \). 

Consid\'erons $\psi:y\mapsto y^E:Y\to X$. Notons que, suvant le 
corollaire de la th\'eor\`eme \ref{thmayet}, les \'el\'ements de la 
base de $X,\C$ ne sont que les $\sigma(e)$, $e\in E=P\cup Q\cup 
\{0,1\}$. Ainsi, il faut consid\'erer les cas suivants: $e=0$ 
ou $e=1$, $e\in P$, et $e\in Q$. On a 
\begin{equation} 
\begin{array}{rcccl}
\psi^{-1}(\sigma(0))&=&
\{y_1\}&=&
\bigcap\{\sigma_1(p)\mid\quad p\in P\}\in C(Y)\cr
\psi^{-1}(\sigma(1))&=&
\{y_0\}&=&
\bigcap\{\sigma_0(p)\mid\quad p\in P\}\in C(Y) \cr
\psi^{-1}(\sigma(p))&=&\sigma_1(p)&&\forall p\in P \cr 
 \psi^{-1}(\sigma(q))&=&\sigma_0(\rho^{-1}(q))&&\forall q\in Q 
\end{array} \label{hom} 
\end{equation} 
et on voit que les es\-pa\-ces de fer\-me\-tu\-re $X,\C$ et $Y,\C$ 
sont en effet ho\-m\'e\-o\-mor\-phi\-ques.  \endproof 

Il reste de souligner les \'el\'ements de $P$ entre les parties 
$CO(Y)$. Rappelons que $Y$ est encore muni de la couple des 
fermetures $\C_1,\C_2$ (\ref{c12}). 

\begin{lemma} L'ensemble ordonn\'e $P$ est isomorhique \`a la 
famille des parties {\em propres} \`a la fois $\C_1$-ferm\'ees et 
$\C_2$-ouvertes de l'espace $Y$:  
\[ (P,\le)\quad\simeq\quad (C_1O_2(Y)\setminus\{\emptyset, 
Y\},\subseteq) \] 
\end{lemma} 

\paragraph{Epreuve.} On montra que cet isomorhisme est 
accompli par $\sigma_1:P\to\exp Y$ et son inverse. Rappelons encore 
que les \'el\'ements de la base de $X,\C$ ne sont que les 
$\sigma(e)$, $e\in E=P\cup Q\cup \{0,1\}$.  Il s'ensuit de 
(\ref{hom}) les preimages des \'el\'ements de $P$ ne sont que les 
$\sigma_1(p)$, $p\in P$. Mais pour tout $p\in P$ la partie 
$Y\setminus\sigma_1(p)=\{y\in Y\mid\quad y(p)=0\}=\sigma_0(p)$, 
alors $\sigma_1(p)\in O_2(X)$.  Pour montrer que c'est 
l'isomorphisme notons que $\sigma_1$ est la restriction de 
l'isomorphisme (Th\'eor\`eme \ref{thmayet}) $\sigma$ sur une partie 
$P$ de $E$. 
\endproof 

\section{Le cas particulier: EOC} 

Etudions la relation de la repr\'esentation introduite ici avec 
telle concernante les EOC. Soit $P$ un EOC avec la 
compl\'ementation $p\mapsto p^c$. Notons que la 
compl\'ementation $(\cdot)^c$ est aussi d\'efinite sur 
l'EO inverse $Q=P^{op}$ comme $q^c=\rho((\rho^{-1}(q))^c)$. 
Introduisons l'EOC $E$ (\ref{eee}) et posons pour tout $x\in X$ 
(\ref{moreoc}) \begin{equation} \pi(x)=1-x(p^c) \label{pie} 
\end{equation} 

\begin{lemma} Pour tout $x\in X$, $\pi(x):E\to B_2$ est un 
\'el\'ement de $X$, c'est \`a dire, $\pi(x)$ est orthocroissante 
par rapport de la compl\'mentation canonique $(\cdot)^\perp$ sur 
$E$. 
\end{lemma} 

\paragraph{Epreuve.} Pour tout $x\in X$, $\pi(x)$ est \'evidemment 
croissante \'etante la composition des applications 
anticroissantes. Pour tout $p\in P$, $p^\perp{}^c = 
\rho((\rho^{-1}(p^\perp)))^c = p^c{}^\perp$, ainsi 
$\pi(x)(p)+\pi(x)(p^\perp) = 1-x(p^c)+1-x((p^\perp)^c) = 
1+1-(x(p^c)+x^((p^c)^\perp)) = 1$ et $\pi(x)$ est en effet {\em 
ortho}croissante.  
\endproof 

\begin{lemma} L'application $\pi:(X,C_1,C_2)\to (X,C_2,C_1)$ des 
\'espaces \`a deux fermetures est bicontinue. \label{lbicont} 
\end{lemma} 

\paragraph{Epreuve.} V\'erifions cela sur les bases des des 
fermetures appropri\'ees. Soit $A=\sigma_1(p)$ pour quelconque 
$p\in P$, alors $\pi^{-1}(A)=\{x\mid\quad 
\pi(x(p))=1\}=\{x\mid\quad x(p^c)=0\}=\sigma_2(p^c)$. \c Ca 
signifie que pour tout $A\in C_1O_2(X)$ nous avons $\pi^{-1}(A)\in 
C_2O_1(X)$. 
\endproof 

L'espace $X=\mor_{EO}(P)$ est plus grand que $Y=\mor_{EOC}(P)$, et 
on obtient la charact\'erization suivante: 

\begin{lemma} $Y$ est l'ensemble des points fix\'ees $\fix(\pi)$ de 
l'application $\pi:X\to X$. \label{lfix} 
\end{lemma} 

\paragraph{Epreuve.} Soit $y\in Y$, alors pour tout $p\in P$ 
on a $\pi(y)(p)=1-y(p^c)=1-(1-y(p))=y(p)$. Conversement, soit 
$y\in\fix(\pi)$, alors pour tout $p\in P$ on a $y(p)+y(p^c)= 
y(p)+1-\pi(y)(p)=1$. \endproof 

Comme nous avons d\'ej\`a mentionn\'e, l'espace $X$ poss\`ede la 
troisi\`eme fermeture $\C=\C_1\lor\C_2$ (\ref{c12}), qu'on pe\^ut 
restreindre \`a $Y\subseteq X$. Un objet de notre int\'er\^et sera 
l'espace de fermeture $(Y,\C)$. 

\begin{lemma} Si $P$ est un EOC, les EOC $P$ et $CO(Y)$ sont 
isomorphiques. 
\end{lemma} 

\paragraph{Epreuve.} Il suffit de v\'erifier que les fermetures 
$\C_Y=\clos(\sigma(P))$ est $\C=\C_1\lor\C_2$ (\ref{c12}) 
co\"\i ncident. Mais c'est \'evident car ils ont les m\^emes bases 
(\c ca s'ensuit de lemma \ref{lbicont}). \endproof 

Alors, tout est pr\^et pour formuler le th\'eor\`eme suivant. 

\begin{th} Soit $P$ un EO avec l'op\'eration involutive 
$(\cdot)^c:P\to P$, et soit $X=\mor_{EO}(P,B_2)$. Ce $P$ est un EOC 
si et seulement si l'application $\pi$ (\ref{pie}) est le hom\'eomorphisme des 
espaces $(X,\C_1,\C_2)$ et $(X,\C_2,\C_1)$. 
\end{th} 

\paragraph{Epreuve.} La n\'ecessit\'e \'etante d\'ej\`a \'etabli 
dans les lemmas pr\'ec\'edents, et il reste de v\'erifier que la 
condition du lemma est suffisante. Soit $p\in P$ et $q=p\lor p^c$, 
alors $\sigma_2(q) \subseteq \sigma_2(p) \cup \sigma_2(p^c)$. Mais 
$\sigma_2(p^c)=\pi(\sigma_1(p))$, alors $\sigma_2(q) \subseteq 
\sigma_1(p) \cup \sigma_2(p) = \emptyset$, ainsi $q=1$ -- le plus 
grand \'el\'ement de $P$. Au m\^eme fa\c con on montre que $p\land 
p^c = 0$, et on v\'erifie que $(\cdot)^c:P\to P$ fournit les 
compl\'ements. \endproof 

\section{Conclusions} 

Dans cet article on a \'etendu la notion de l'espace dual pour un 
ensemble ordonn\'e compl\'ement\'e dans le cas le plus g\'en\'eral 
des ensembles ordonn\'es arbitraires. Les techniques introduites 
ici peuvent \^etre consid\'er\'ees comme la representation dans un 
sense "naturelle" des EO: on les repr\'esente avec les parties de 
l'universe appropri\'e. De plus, la condition pour une involution 
sur un EO d'\^etre compl\'ementation est introduite ici. 

On exprime la reconnaissance \`a Georges Chevalier et Ren\'e Mayet 
pour les discussions.

\noindent D\'epartment des Math\'ematiques, SPb UEF,\newline 
\noindent Griboyedova 30/32, \newline 
\noindent 191023, St.Petersbourg, Russie 
\end{document}